\documentclass[12pt,reqno]{amsart}

\usepackage{amssymb}
\usepackage[hyperfootnotes=false]{hyperref}

\renewcommand{\phi}{\varphi}

\newcommand{\<}{\langle}
\renewcommand{\>}{\rangle}

\newcommand{\stm}{\setminus}
\newcommand{\seq}{\subseteq}

\newcommand{\longc}{,\ldots,}

\newtheorem{theorem}{Theorem}

\newtheorem{primetheorem}{Theorem}

\newcommand{\refs}[1]{~\ref{s:#1}}
\newcommand{\reft}[1]{~\ref{t:#1}}
\newcommand{\refe}[1]{~\eqref{e:#1}}

\title[Small asymmetric sumsets]%
  {Small asymmetric sumsets \\ in elementary abelian $2$-groups}

\author{Chaim Even Zohar}
\address{Einstein Institute of Mathematics, The Hebrew University,
        Jerusalem 91904, Israel}
\email{chaim.evenzohar@mail.huji.ac.il}

\author{Vsevolod F. Lev}
\address{Department of Mathematics, The University of Haifa at Oranim,
  Tivon 36006, Israel}
\email{seva@math.haifa.ac.il}

\begin{document}
\baselineskip=16pt

\maketitle

\section{Summary}
Let $A$ and $B$ be subsets of an elementary abelian $2$-group $G$, none of
which are contained in a coset of a proper subgroup. Extending onto
potentially distinct summands a result of Hennecart and Plagne, we show that
if $|A+B|<|A|+|B|$, then either $A+B=G$, or the complement of $A+B$ in $G$ is
contained in a coset of a subgroup of index at least $8$ (whence
$|A+B|\ge\frac78\,|G|$). We indicate conditions for the containment to be
strict, and establish a refinement in the case where the sizes of $A$ and $B$
differ significantly.
\renewcommand{\thefootnote}{}
\footnote{2010 \emph{Mathematics Subject Classification}: Primary 11P70; Secondary 05D99.}
\footnote{\emph{Key words and phrases}: sumset, abelian 2-group, affine span.}
\renewcommand{\thefootnote}{\arabic{footnote}}
\setcounter{footnote}{0}

\section{Background and introduction}\label{s:introduction}

For subsets $A$ and $B$ of an abelian group, we denote by $A+B$ the sumset of
$A$ and $B$:
  $$ A+B := \{ a+b\colon a\in A,\ b\in B \}. $$
We abbreviate $A+A$ as $2A$. By $\<A\>$ we denote the affine span of $A$
(which is the smallest coset that contains $A$).

Pairs of finite subsets $A$ and $B$ of an abelian group with $|A+B|<|A|+|B|$
are classified by the classical results of Kneser and Kemperman
\cite{kneser,kemperman}. Recursive in its nature, this classification is
rather complicated in general, but it has been observed that for the special
case where the underlying group is an elementary abelian $2$-group (that is,
a finite abelian group of exponent $2$), explicit closed-from results can be
obtained. Particularly important in our present context is the following
theorem due to Hennecart and Plagne.
\begin{theorem}[{\cite[Theorem 1]{hennecart_plagne}}]\label{t:hp}
Let $A$ be a subset of an elementary abelian $2$-group $G$ such that
$\<A\>=G$. If $|2A|<2|A|$, then either $2A=G$, or the complement of $2A$ in
$G$ is a coset of a subgroup of index at least $8$. Consequently,
$|2A|\ge\frac78\,|G|$.
\end{theorem}

We mention two directions in which Theorem \reft{hp} was later developed.
First, in connection with Freiman's structure theorem, much attention has
been attracted to the function $F$ defined by
  $$ F(K) := \sup \{ |\<A\>|/|A|\colon |2A|\le K|A| \},\quad K\ge 1 $$
where $A$ runs over non-empty subsets of elementary abelian $2$-groups. It is
not difficult to derive from Theorem \reft{hp} that
  $$ F(K) = \begin{cases}
               K &\text{if $1\le K<\frac74$} \\
               \frac87\,K  &\text{if $\frac74\le K<2$}
            \end{cases}; $$
this is, essentially, \cite[Corollary~2]{hennecart_plagne}.  A result of
Ruzsa \cite{ruzsa_thm} shows that $F(K)$ is finite for each $K\ge 1$ and
indeed, $F(K)\le K^22^{K^4}$. Various improvements for $K\ge 2$ were obtained
by Deshouillers, Hennecart, and Plagne \cite{dhp}, Sanders \cite{sanders},
Green and Tao \cite{green_tao}, and Konyagin \cite{konyagin}, and the exact
value of $F(K)$ was eventually established in \cite{evenzohar}.

In another direction,  \cite[Theorem~5]{lev_structure} establishes the
precise structure of those subsets $A$ satisfying $|2A|<2|A|$ --- in contrast
with Theorem \reft{hp} which describes the structure of the sumset $2A$ only.

The goal of the present paper is to extend Theorem \reft{hp} onto addition of
two potentially distinct set summands. In this case the assumption
$|A+B|<|A|+|B|$ does not guarantee any longer that the complement of $A+B$ is
a coset of a subgroup of index at least $8$, as evidenced, for instance, by
the following construction: represent the underlying group $G$ as a direct
sum $G=H\oplus F$ with $|H|=8$, fix a generating set
 $\{h_1,h_2,h_3\}\subset H$ and an arbitrary proper subset $F_0\subsetneq F$,
and let
\begin{align*}
  A &:= \big( \{h_1,h_2,h_3\}+F \big) \cup \{0\}, \\
  B &:= \big( \{h_1+h_2,h_2+h_3,h_3+h_1,h_1+h_2+h_3\}+F \big) \cup F_0.
\end{align*}
The complement of $A+B$ in $G$ is easily verified to be the complement of
$F_0$ in $F$, which need not be a coset, and
  $$ |A+B|=|G|-(|F|-|F_0|) = |A| + |B| - 1. $$
It turns out, however, that while the complement of $A+B$ may fail to be a
coset of a subgroup of index at least $8$, it is necessarily \emph{contained}
in a such a coset --- and indeed, in a coset of a subgroup of larger index if
the summands differ significantly in size.

For subsets $A$ and $B$ of an abelian group and a group element $g$, let
$\nu_{A,B}(g)$ denote the number or representations of $g$ in the form
$g=a+b$ with $a\in A$ and $b\in B$, and let
  $$ \mu_{A,B}:=\min\{\nu_{A,B}(g)\colon g\in A+B \}. $$
The following theorem, proved in Section \refs{proof}, is our main result.
\begin{theorem}\label{t:main}
Let $A$ and $B$ be subsets of an elementary abelian $2$-group $G$ such
that $\<A\>=\<B\>=G$. If $|A+B|<\min\{|A|+|B|,|G|\}$, then the complement
of $A+B$ in $G$ is contained in a coset of a subgroup of index $8$.
Moreover, if $\mu_{A,B}=1$, then the containment is strict.
\end{theorem}

We could get a stronger conclusion in the ``highly asymmetric" case.
\begin{primetheorem}\label{t:asymmetric}
Let $A$ and $B$ be subsets of an elementary abelian $2$-group $G$ such
that $\<A\>=\<B\>=G$. If $|A+B|<\min\{|A|+|B|,|G|\}$ and
$|B|\ge\Big(1-\frac{k+1}{2^k}\Big)\,|G|$ with integer $k\ge 4$, then the
complement of $A+B$ in $G$ is contained in a coset of a subgroup of index
$2^k$. Moreover, if $\mu_{A,B}=1$, then the containment is strict.
\end{primetheorem}

Notice that in the statements of Theorems \reft{main} and
\reft{asymmetric} we disposed of the case where the sumset $A+B$ is the
whole group by assuming from the very beginning that $|A+B|<|G|$.

The bounds on the subgroup index in Theorems \reft{main} and
\reft{asymmetric} are best possible under the stated assumptions. To see
this, fix an integer $k\ge 3$ (the case $k=3$ addressing Theorem
\reft{main}), consider a decomposition $G=H\oplus F$ with $|H|=2^k$, choose a
generating set $\{0,h_1\longc h_k\}\subset H$ and two arbitrary elements
$g_1,g_2\in G$, and let
\begin{align*}
  A &:= g_1 + \{0,h_1\longc h_k\} + F, \\
  B &:= g_2 + ( H\stm \{0,h_1\longc h_k\} ) + F.
\end{align*}
Then $|B|=\Big(1-\frac{k+1}{2^k}\Big)|G|$, the complement of $A+B$ in $G$ is
$g_1+g_2+F$, and
  $$ |A+B| = |G|-|F| = |A|+|B|-|F|. $$
Indeed, analyzing carefully the argument in Section \refs{proof}, one can see
that if $B$ is not of the form just described, then the containment in the
conclusion of Theorem \reft{asymmetric} is strict.

An almost immediate corollary of Theorem \reft{main} is that if $A$ and
$B$ are subsets of an elementary abelian $2$-group $G$ such that
$\<A\>=\<B\>=G$ and $|A+B|<\frac78\,(|A|+|B|)$, then $A+B=G$. In fact,
Kneser's theorem \cite{kneser} yields a stronger result: if
$\<A\>=\<B\>=G$ and $|A+B|<|A|+\frac34\,|B|$, then $A+B=G$. Omitting the
proof, which is nothing more than a routine application of Kneser's
theorem, we confine ourselves to the remark that both assumptions
$\<A\>=G$ and $\<B\>=G$ are crucial. This follows by considering the
situation where $B$ is an index-$8$ subgroup of $G$, and $A$ is a union
of $4$ cosets of $B$ (which is not a coset itself), and that where $A$ is
an index-$4$ subgroup, and $B$ is a union of three cosets of $A$.

We deduce Theorems \reft{main} and \reft{asymmetric} from
\cite[Theorem~2]{lev_structure}, quoted in the next section as
Theorem~\reft{lev}. Based on the well-known Kemperman's structure
theorem, this result establishes the structure of pairs $(A,B)$ of
subsets of an abelian group such that $|A+B|<|A|+|B|$. The deduction of
Theorems \reft{main} and \reft{asymmetric} from Theorem~\reft{lev} is
presented in Section \refs{proof}.

\section{Pairs of sets with a small sumset}\label{s:toolbox}

The contents of this section originate from \cite{kemperman} and
\cite{lev_structure}. Our goal here is to introduce
\cite[Theorem~2]{lev_structure}, from which Theorems \reft{main} and
\reft{asymmetric} will be derived in the next section.

For a subset $A$ of the abelian group $G$, the (maximal) period of $A$ will
be denoted by $\pi(A)$; recall that this is the subgroup of $G$ defined by
  $$ \pi(A) := \{g\in G\colon A+g=A \}, $$
and that $A$ is called \emph{periodic} if $\pi(A)\ne\{0\}$ and
\emph{aperiodic} otherwise.

By an arithmetic progression in the abelian group $G$ with difference
 $d\in G$, we mean a set of the form $\{g+d,g+2d\longc g+nd\}$, where $n$ is a
positive integer.

Essentially following  Kemperman's paper \cite{kemperman}, we say that the
pair $(A,B)$ of finite subsets of the abelian group $G$ is \emph{elementary}
if at least one of the following conditions holds:
\begin{itemize}
\item[(I)]   $\min\{|A|,|B|\}=1$;
\item[(II)]  $A$ and $B$ are arithmetic progressions sharing a common
    difference, the order of which in $G$ is at least $|A|+|B|-1$;
\item[(III)] $A=g_1+(H_1\cup\{0\})$ and $B=g_2-(H_2\cup\{0\})$, where
    $g_1,g_2\in G$, and where $H_1$ and $H_2$ are non-empty subsets of a
    subgroup $H\le G$ such that $H=H_1\cup H_2\cup\{0\}$ is a partition
    of $H$; moreover, $c:=g_1+g_2$ is the unique element of $A+B$ with
    $\nu_{A,B}(c)=1$;
\item[(IV)] $A=g_1+H_1$ and $B=g_2-H_2$, where $g_1,g_2\in G$, and where
    $H_1$ and $H_2$ are non-empty, aperiodic subsets of a subgroup $H\le
    G$ such that $H=H_1\cup H_2$ is a partition of $H$; moreover,
    $\mu_{A,B}\ge 2$.
\end{itemize}

Notice, that for elementary pairs of type (III) we have $|A|+|B|=|H|+1$,
whence $A+B=g_1+g_2+H$ by the box principle. Also, for type (IV) pairs we
have $|A|+|B|=|H|$ and $A+B=g_1+g_2+(H\stm\{0\})$; the reader can consider
the latter assertion as an exercise or find a proof in \cite{lev_structure}.

We say that the pair $(A,B)$ of subsets of an abelian group satisfies
\emph{Kemperman's condition} if
\begin{equation}\label{e:kemp-cond}
  \text{either}\ \pi(A+B)=\{0\},\ \text{or}\ \mu_{A,B}=1.
\end{equation}

Given a subgroup $H$ of the abelian group $G$, by $\phi_H$ we denote the
canonical homomorphism from $G$ onto the quotient group $G/H$.

We are at last ready to present our main tool.
\begin{theorem}[{\cite[Theorem~2]{lev_structure}}]\label{t:lev}
Let $A$ and $B$ be finite, non-empty subsets of the abelian group $G$. A
necessary and sufficient condition for $(A,B)$ to satisfy both
  $$ |A+B| < |A| + |B| $$
and Kemperman's condition \refe{kemp-cond} is that either $(A,B)$ is an
elementary pair, or there exist non-empty subsets $A_0\seq A$ and $B_0\seq B$
and a finite, non-zero, proper subgroup $F<G$ such that
\begin{itemize}
\item[(i)]   each of $A_0$ and $B_0$ is contained in an $F$-coset,
    $|A_0+B_0|=|A_0|+|B_0|-1$, and the pair $(A_0,B_0)$ satisfies
    Kemperman's condition;
\item[(ii)]  each of $A\stm A_0$ and $B\stm B_0$ is a (possibly empty)
    union of $F$-cosets;
\item[(iii)] the pair $(\phi_F(A),\phi_F(B))$ is elementary; moreover,
    $\phi_F(A_0)+\phi_F(B_0)$ has a unique representation as a sum of an
    element of $\phi_F(A)$ and an element of $\phi_F(B)$.
\end{itemize}
\end{theorem}

\section{Proof of Theorems \reft{main} and \reft{asymmetric}}\label{s:proof}

We give Theorems \reft{main} and \reft{asymmetric} one common proof.

If $|G|\le 4$, then the assumption $\<A\>=\<B\>=G$ implies $A+B=G$, and
we therefore assume $|G|\ge 8$ and use induction on $|G|$.

If Kemperman's condition \refe{kemp-cond} fails to hold, then, in particular,
$H:=\pi(A+B)$ is a non-zero subgroup. In this case we observe that the
assumptions $\<A\>=\<B\>=G$ and $|A+B|<|G|$ imply
$\<\phi_H(A)\>=\<\phi_H(B)\>=G/H$ and $|\phi_H(A)+\phi_H(B)|<|G/H|$,
respectively, and
\begin{equation}\label{e:Blarge}
  |B|\ge\Big(1-\frac{k+1}{2^k}\Big)\,|G|
\end{equation}
implies $|\phi_H(B)|\ge\Big(1-\frac{k+1}{2^k}\Big)\,|G/H|$. Hence, by the
induction hypothesis, the complement of $\phi_H(A)+\phi_H(B)=\phi_H(A+B)$ in
$G/H$ is contained in a coset of a subgroup of index $8$ and indeed, of index
$2^k$ under the assumption \refe{Blarge}, and so is the complement of $A+B$
in $G$.

From now on we assume that Kemperman's condition \refe{kemp-cond} holds true,
and hence Theorem~\reft{lev} applies.

If $(A,B)$ is an elementary pair in $G$, then it is of type III or IV, in
view of the assumptions $|G|\ge 8$ and $\<A\>=\<B\>=G$. Moreover, by the same
reason, the subgroup $H\le G$ in the definition of elementary pairs is, in
fact, the whole group $G$.  We conclude that $(A,B)$ is actually of type IV:
for, if it were of type III, we would have $A+B=G$ (see a remark after the
definition of elementary pairs). Consequently, $\mu_{A,B}\ge 2$ and the
complement of $A+B$ in $G$ is a singleton; that is, a coset of the zero
subgroup. To complete the treatment of the present case, we denote by $n$ the
rank of $G$ and notice that \refe{Blarge} implies $|A|=|G|-|B|\le
(k+1)2^{n-k}$, while $\<A\>=G $ gives $|A|\ge n+1$. Hence,
$(n+1)/2^n\le(k+1)/2^k$. As a result, $n\ge k$, and therefore the zero
subgroup has index $|G|\ge2^k$.

Finally, consider the situation where $(A,B)$ is not an elementary pair in
$G$, and find then $A_0\seq A,\, B_0\seq B$, and $F<G$ as in the conclusion
of Theorem~\reft{lev}. Observe that $\<\phi_F(A)\>=\<\phi_F(B)\>=G/F$ yields
$\min\{|\phi_F(A)|,|\phi_F(B)|\}\ge 2$, so that $(\phi_F(A),\phi_F(B))$
cannot be an elementary pair in $G/F$ of type I or II. Indeed,
$(\phi_F(A),\phi_F(B))$ cannot be of type IV either, as in this case we would
have $\mu_{\phi_F(A),\phi_F(B)}\ge 2$, contrary to Theorem~\reft{lev}~(iii).
Thus, $(\phi_F(A),\phi_F(B))$ is of type III, and
$\<\phi_F(A)\>=\<\phi_F(B)\>=G/F$ implies that the subgroup of the quotient
group $G/F$ in the definition of elementary pairs is actually the whole group
$G/F$. As a result, we derive from Theorem \reft{lev} that the complement of
$A+B$ in $G$ is the complement of $A_0+B_0$ in the appropriate $F$-coset.

Write $|G/F|=2^m$; to complete the proof it remains to show that $m\ge 3$,
and if \refe{Blarge} holds then, indeed, $m\ge k$. To this end we notice that
$\<\phi_F(A)\>=\<\phi_F(B)\>=G/F$ gives $\min\{|\phi_F(A)|,|\phi_F(B)|\}\ge
m+1$; compared to $|\phi_F(A)|+|\phi_F(B)|=2^m+1$, this results in $2m+2\le
2^m+1$, whence $m\ge 3$. Finally, $|\phi_F(B)|\ge(1-(k+1)/2^k)\,2^m$ gives
$|\phi_F(A)|\le(k+1)2^{m-k}+1$. Combined with $|\phi_F(A)|\ge m+1$ this leads
to $m\le (k+1)2^{m-k}$. As the right-hand side is a decreasing function of
$k$, if we had $m<k$, the last inequality would yield $m\le
(m+2)2^{m-(m+1)}$, which is wrong.

Note that the condition $\mu_{A,B}=1$ can hold only under the last
scenario (where $(A,B)$ is not an elementary pair in $G$). As we have
shown, in this case the complement of $A+B$ is strictly contained in an
$F$-coset, and the strict containment assertion follows. \qed

\section*{Acknowledgment}

The first author would like to thank his thesis advisor, Professor Nati
Linial, for a patient and helpful guidance.

\bibliographystyle{amsalpha}
\bibliography{SmallF2n}

\end{document}